\documentclass[reqno,12pt]{amsart}
\usepackage{amsfonts}
\usepackage{amssymb}
\usepackage{bbm}
\usepackage{times}

\def\D{\mathbb{D}}

\def\C{\mathbb{C}}

\def\N{\mathbb N}

\def\msk{\medskip}

\def\bege{\begin{equation}} \def\ende{\end{equation}}

\def\begr{\begin{eqnarray}} \def\endr{\end{eqnarray}}

\def\bege{\begin{equation}} \def\ende{\end{equation}}
\def\begr{\begin{eqnarray}} \def\endr{\end{eqnarray}}
\def\bnum{\begin{enumerate}} \def\enum{\end{enumerate}}

\begin{document}

\title[Differences of composition operators]{Differences of composition operators between different Hardy spaces}
\author{Yecheng Shi and Songxiao Li$^*$ }

\address{Yecheng Shi, Faculty of Information Technology,
Macau University of Science and Technology, Avenida Wai Long,
Taipa, Macau.
}\email{ 09ycshi@sina.cn}

\address{Songxiao Li\\ Institute of Fundamental and Frontier Sciences, University of Electronic Science and Technology of China,
610054, Chengdu, Sichuan, P.R. China.  \newline
Institute of Systems Engineering, Macau University of Science and Technology, Avenida Wai Long, Taipa, Macau. } \email{jyulsx@163.com}

\begin{abstract}  In this paper, we give some estimates for the norm and essential norm of the differences of two composition operators between different Hardy spaces.
\thanks{*Corresponding author.}
\end{abstract}

\maketitle

\section{Introduction}

Let $\D$ denote the open unit disk of the complex plane $\C$. We denote the closure and the unit circle of $\D$ by $\overline\D$ and $\partial\D$, respectively. Let $H(\D)$ be  the class of functions analytic in $\D$. Let $dm=\frac{d\theta}{2\pi}$  denote the normalized Lebesgue measure on $\partial\D$.
The Lebesgue space $L^p(m)$ will also be denoted by $L^p(\partial\D)$, $0<p<\infty$.
For $0<p<\infty,$ let $H^p$ denote the Hardy space of all  $f\in H(\D)$ such that
$$\|f\|_{p}^p=\sup_{0<r<1} \int_{ \partial \D}|f(r\xi)|^pdm(\xi)<\infty.$$

Recall that if $f\in H^p(\D)$, then the radial limits $\lim_{r\to1}f(re^{i\theta})$ exist almost everywhere on $\partial\D$ and will be denoted also by $f$, which belongs to $L^p(\partial\D)$  and
$$\|f\|_{p}^p=\frac{1}{2\pi}\int_{0}^{2\pi}|f(e^{i\theta})|^pd\theta.$$
The space $H^\infty(\D)$ consists of all bounded analytic functions on $\D$, and its norm is given by the supremum norm on $\D$.

For $a\in\D$, let $\sigma_{a}(z):=\frac{a-z}{1-\bar{a}z}$ be the disc automorphism that exchanging $0$ for $a$. Let $\bigtriangleup(a,r) := \{z \in\D: |\sigma_a(z)| < r\}$  denote the pseudohyperbolic disk centered at $a$ with radius $r$.
For two points $z$, $w\in\D$, the pseudohyperbolic distance is given by
$$\rho(z,w)=|\sigma_{w}(z)|=\big|\frac{z-w}{1-\bar{w}z}\big|.$$

Let $\varphi$ and $\psi$ be two analytic self-maps of $\D$. We write that $\sigma(z)=\rho(\varphi(z),\psi(z)).$
Note that $\sigma $ also has a radial extension $\sigma^{*}$ almost everywhere on $\partial\D$. Indeed, if $\varphi\neq\psi$, then the radial limits of $\varphi$ and $\psi$ can  coincide only on a set of measure zero.  We will use the same notation for the function $\sigma$ and its radial extension.

Let $\varphi$ be an analytic self-map of $\D$.  The composition operator $C_\varphi$ on $H(\D)$ is defined by
$C_\varphi f=f\circ\varphi$. Furthermore, if $u$ is a Borel measurable function, a weighted composition
operator $uC_{\varphi}$ on $H(\D)$ is defined by
$$
(uC_\varphi f)(z) = u(z)f( \varphi(z)), ~~~~ ~f \in H(\D), ~\\~z\in \D.
$$

The study of the differences of  composition operators was
started on the Hardy space $H^2$. The main purpose for this study   is to
understand the topological structure of the set of composition
operators $\mathcal{C}(H^2)$, see \cite{Be,cm,Sj}.  After that, such
related problems have been studied on several spaces of analytic
functions by many authors, see for example \cite{GGNS, li12,LS,MOZ, Mo,Ni,NS,S1,S2,SL}.
Motivated by \cite{Mo,S1,S2}, in this paper, we study the differences of composition operators between different Hardy spaces.

In the Hardy spaces setting, Goebeler \cite{G} showed that for $0< q<p<\infty$, $C_\varphi-C_\psi:H^p\to H^q$ is compact if and only if the composition operators $C_\varphi$ and $C_\psi$ between these spaces are both compact, that is,
$|\varphi|<1$ and $|\psi|<1$ a.e. on $\partial\D$.
Nieminen and Saksman \cite{NS} proved that $C_\varphi-C_\psi$ is compact on $H^p$ for all $  p \in [1,\infty)$ if and only if $C_\varphi-C_\psi$ is compact for some $  p \in [1,\infty)$.
But the complete characterization of the compactness of $C_\varphi-C_\psi:H^p\to H^p$ is still open.

Recently, Saukko in \cite{S1} asked the following question: {\it is boundedness and compactness of the difference operator in Hardy spaces enough to guarantee the boundedness and compactness of corresponding weighted composition operators?}  In this paper, we give partly
positive answers to this question. The main results of this paper are the following theorems.\msk

\noindent{\bf Theorem 1.1.}  {\it Let $1<p<q< \infty$. Suppose $\varphi$ and $\psi$ are analytic self-maps of $\D$. Then $C_\varphi-C_\psi:H^p\to H^q$ is bounded if and only if both weighted composition operators $\sigma C_\varphi$ and $\sigma C_\psi$ map $H^p$ into $H^q$. Furthermore,
\begr
&&\|C_\varphi-C_\psi\|_{H^p\to H^q}^q\nonumber\\
&\approx&\sup_{a\in\D}\int_{\partial{\D}} \Big|\bigg(\frac{1-|a|^2}{(1-\overline{a}\varphi(\xi))^2}\bigg)^{\frac{1}{p}}-
\bigg(\frac{1-|a|^2}{(1-\overline{a}\psi(\xi))^2}\bigg)^{\frac{1}{p}}\Big|^qdm(\xi)\nonumber\\
&\approx&\sup_{a\in\D}\int_{\partial{\D}} \Big(\frac{1-|a|^2}{|1-\overline{a}\varphi(\xi)|^2}+\frac{1-|a|^2}{|1-\overline{a}\psi(\xi)|^2}\Big)^{\frac{q}{p}}\rho(\varphi(\xi),\psi(\xi))^qdm(\xi).
\nonumber\endr}

\noindent{\bf Theorem 1.2.}  {\it Let $1<p<q< \infty$. Suppose $\varphi$ and $\psi$ are analytic self-maps of $\D$ such that $C_\varphi-C_\psi: H^p\to H^q$ is bounded. Then
\begr
&&\|C_\varphi-C_\psi\|_{e, H^p\to H^q}^q\nonumber\\
&\approx&\limsup_{|a|\to 1}\int_{\partial{\D}} \Big|\bigg(\frac{1-|a|^2}{(1-\overline{a}\varphi(\xi))^2}\bigg)^{\frac{1}{p}}-
\bigg(\frac{1-|a|^2}{(1-\overline{a}\psi(\xi))^2}\bigg)^{\frac{1}{p}}\Big|^qdm(\xi)\nonumber\\
&\approx&\limsup_{|a|\to 1}\int_{\partial{\D}} \Big(\frac{1-|a|^2}{|1-\overline{a}\varphi(\xi)|^2}+\frac{1-|a|^2}{|1-\overline{a}\psi(\xi)|^2}\Big)^{\frac{q}{p}}
\rho(\varphi(\xi),\psi(\xi))^qdm(\xi).\nonumber\endr}

Recall that the essential norm of a bounded linear operator $T:X\rightarrow Y$ is its distance to the set of
compact operators $K$ mapping $X$ into $Y$, that is,
$$\|T\|_{e, X\rightarrow Y}=\inf\{\|T-K\|_{X\rightarrow Y}: K~\mbox{is compact}~~\},$$ where  $X,Y$ are Banach spaces and
$\|\cdot\|_{X\rightarrow Y}$ is the operator norm.

The present paper is organized as follows. In Section 2,  we study  weighted composition operators between Hardy spaces.
 In Section 3,  we state some lemmas and give the proofs of Theorems 1.1 and 1.2.

For two quantities $A$ and $B$, we use the abbreviation $A\lesssim B$ whenever there is a positive constant $c$ (independent of the associated variables) such that $A\leq cB$. We write $A\approx B$, if $A\lesssim B\lesssim A$.

\section{weighted composition operators from $H^p$ to $H^q$}

In this section, we collect some characterizations of weighted composition operators between different Hardy spaces.
Given any measure $\mu$ on $\overline\D$, we denote by $\mu|_{\D}$ and $\mu|_{\partial\D}$ its restrictions to the Borel subsets of $\D$ and $\partial\D$, respectively. By Lemma 2.1 of \cite{BJ}, the $s$-Carleson measure on $\overline\D$ is defined as follows.\msk

\noindent{\bf Definition 2.1.}  {\it Let $0<p, q<\infty$, $\mu$ be a Borel measure on $\overline\D$. Then the measure $\mu$ is call a $\frac{q}{p}$-Carleson measure on $\overline\D$ if the inclusion map $I_\mu:H^p\to L^q(\mu,\overline\D)$ is bounded, i.e.,  there exists a constant $C$ such
that
\begr
\Big(\int_{\overline\D}|f(z)|^qd\mu(z)\Big)^{\frac{1}{q}}\leq C\|f\|_{p}\nonumber
\endr
for every $f\in H^p$. Furthermore, $\mu$ is a vanishing $\frac{q}{p}$-Carleson measure on $\D$ if the inclusion map
$I_{\mu}: H^p\to L^q(\mu,\overline\D)$ is compact.
}\msk

For an interval $I\subset \partial\D$,  the Carleson square is defined by
$$S(I)=\{re^{it}\in\D:1-|I|<r<1, e^{it}\in I\},$$
 where $|E|$ denotes the Lebesgue measure of the measurable set $E\subset\partial\D$. If $a\in \D\backslash\{0\}$, let  $I_a=\{e^{i\theta}:|arg(ae^{-i\theta})|\leq \frac{1-|a|}{2}\}$, and denote
$S(a)=S(I_a)$. For convenience, we put $I_0=\partial\D$ and $S(0)=\D$.

Suppose $u:\partial\D\to\C$ is a measurable function and $\varphi$ is an analytic self-map of $\D$. Define the measure $\mu_{u,\varphi}$ in $\overline\D$ by
$$\mu_{u,\varphi}(E)=\int_{\varphi^{-1}(E)\cap\partial\D}|u(z)|^qdm(z)  $$
for all Borel set $E\subset\overline\D$.

We   need the following results about weighted composition operators on Hardy spaces from \cite{BJ} and \cite{CZ2}.\msk

\noindent{\bf Theorem 2.1.}  {\it Suppose $0<p<q< \infty$ and $0<r<1$. Let $u:\partial\D\to\C$ be a measurable function and $\varphi$ an analytic self-map of $\D$. Then the following statements are equivalent:

(i) The weighted composition operator $uC_\varphi:H^p\to L^q(\partial\D)$ is bounded.

(ii) $\mu_{u,\varphi}|_{\partial\D}=0$ and $\|\mu_{u,\varphi}|_{\D}\|_{p,q}^q:=\sup_{a\in\D}\frac{\mu(S(a))}{(1-|a|^2)^{\frac{q}{p}}}<\infty$.

(iii) $\mu_{u,\varphi}|_{\partial\D}=0$ and $\|\mu_{u,\varphi}|_{\D}\|_{p,q,r}^q:=\sup_{a\in\D}\frac{\mu(\bigtriangleup(a,r))}{(1-|a|^2)^{\frac{q}{p}}}<\infty$.

(iv) $\sup_{a\in\D}\int_{\overline\D}\big(\frac{1-|a|^2}{|1-\overline{a}z|^2}\big)^\frac{q}{p}d\mu_{u,\varphi}(z)<\infty.$

(v) $\sup_{a\in\D}\|(uC_\varphi)k_a\|_{L^q(\partial\D)}^q<\infty,$ where $k_a(z)=\big(\frac{1-|a|^2}{(1-\overline{a}z)^2}\big)^{\frac{1}{p}}$.

Furthermore, $\|uC_\varphi\|_{H^p\to L^q(\partial\D)}^q$ and the quantities in (ii), (iii), (iv) and (v) are all comparable
with comparability constants depending only on $p$, $q$  and $r$. }\msk

{\it  Proof.}  The equivalence between $(ii)$ and $(iii)$ can be found in \cite{Lu}.  By the   change of variables, for all $f\in H(\D)$ (see \cite[Lemma 2.1]{CH}, also see \cite{S1}),
 $$\|uC_{\varphi}f\|_{L^q(\partial\D)}=\|f\|_{L^q(\mu_{u,\varphi},\overline\D)}.$$
  Therefore, $uC_\varphi:H^p\to L^q(\partial\D)$ is bounded if and only if the inclusion map $I_{\mu_{u,\varphi}}:H^p\to L^q(\mu_{u,\varphi},\overline\D)$
is bounded, and $\|uC_\varphi\|=\|I_{\mu_{u,\varphi}}\|$. Taking $f=k_a$, we obtain that $(iv)$ and $(v)$ are equivalent.
The equivalence of $(i), (ii)$ and $(iv)$ follows from \cite[Theorem 2.5]{BJ} and the proof of \cite[ Proposition 2.3]{BJ}.
 The comparability of the quantities are clear. The proof is complete. \msk

Let $n\in\N$. Define the partial sum operator $S_n:H(\D)\to H(\D)$ by
$$S_n\Big(\sum_{k=0}^{\infty}a_kz^k\Big)=\sum_{k=0}^na_kz^k.$$
Denote $R_n=I-S_n$. 
For $0<s<1$, we denote $\D_s=\{z\in\D:|z|<s\}$.\msk

\noindent{\bf Theorem 2.2.}  {\it Suppose $1< p<q< \infty$ and $0<r<1$. Suppose that $u:\partial\D\to\C$ is a measurable function and $\varphi$ is an analytic self-map of $\D$ such that the operator $uC_\varphi:H^p\to L^q(\partial\D)$ is bounded. Then

(i) \begr  \|uC_\varphi\|_{e, H^p\to L^q(\partial\D)}^q&\approx&\lim_{s\to 1}\|\mu_{u,\varphi}|_{\D\backslash\D_s}\|_{p,q,r}^q\nonumber\\
 &\approx&\liminf_{n\to\infty}\|(uC_\varphi)R_n\|_{H^p\to L^q(\partial\D)}^q\nonumber\\
 &\approx&\limsup_{n\to\infty}\|(uC_\varphi)R_n\|_{H^p\to L^q(\partial\D)}^q\nonumber\\
 &\approx& \limsup_{|a|\to1}\|(uC_\varphi)k_a\|_{L^q(\partial\D)}^q\nonumber\\
   &\approx&\limsup_{|a|\to1}\frac{\mu_{u,\varphi}(\bigtriangleup(a,r))}{(1-|a|^2)^{\frac{q}{p}}}.  \nonumber
 \endr

(ii) For every $0<\eta<1$,
$$\lim_{n\to\infty}\sup_{\|f\|_{p}\leq1}\int_{\varphi^{-1}(\D_\eta)}|(C_\varphi\circ R_nf)(\xi)|^qdm(\xi)=0.$$ }

{\it Proof.} (i) First, we prove that
$$\lim_{s\to1}\|\mu_{u,\varphi}|_{\D\backslash\D_s}\|_{p,q,r}^q\approx
\limsup_{|a|\to1}\frac{\mu_{u,\varphi}(\bigtriangleup(a,r))}{(1-|a|^2)^{\frac{q}{p}}}.$$
Let $$t_r(s)=\frac{s-r}{1-sr}.$$
 After a calculation, we get that $\bigtriangleup(a,r)\cap(\D\backslash\D_s)\neq0$ if and only if $|a|\geq t_r(s)$.
 It is easy to see that $t_r(s)$ is continuous and increasing on $[r,1)$, and $\lim_{s\to1}t_r(s)=1$.
Thus,
\begr
\limsup_{|a|\to1}\frac{\mu_{u,\varphi}(\bigtriangleup(a,r))}{(1-|a|^2)^{\frac{q}{p}}}&=&\lim_{s\to1}\sup_{|a|\geq t_r(s)}\frac{\mu_{u,\varphi}(\bigtriangleup(a,r))}{(1-|a|^2)^{\frac{q}{p}}}\nonumber\\
&\geq&\lim_{s\to1}\sup_{|a|\geq t_r(s)}\frac{\mu_{u,\varphi}(\bigtriangleup(a,r)\cap(\D\backslash\D_s))}{(1-|a|^2)^{\frac{q}{p}}}\nonumber\\
&=&\lim_{s\to1}\sup_{a\in\D}\frac{\mu_{u,\varphi}(\bigtriangleup(a,r)\cap(\D\backslash\D_s))}{(1-|a|^2)^{\frac{q}{p}}}\nonumber\\
&=&\lim_{s\to1}\|\mu_{u,\varphi}|_{\D\backslash\D_s}\|_{p,q,r}^q.\nonumber
\endr
Denote $A=\lim_{s\to1}\|\mu_{u,\varphi}|_{\D\backslash\D_s}\|_{p,q,r}^q$. For any $\epsilon>0$, there exists $0<t<1$, such that if $t\leq s<1$, we have
$$\|\mu_{u,\varphi}|_{\D\backslash\D_s}\|_{p,q,r}^q<A+\epsilon.$$
 For any fixed $s$ ($0<s<1$), we know that $\bigtriangleup(a,r)\subset\D\backslash\D_s$, as $|a|$ close enough to 1. Therefore, there exists a $l$, $0<l<1$, such that
\begr
\|\mu_{u,\varphi}|_{\D\backslash\D_s}\|_{p,q,r}^q&=&\sup_{a\in\D}\frac{\mu_{u,\varphi}(\bigtriangleup(a,r)\cap(\D\backslash\D_s))}{(1-|a|^2)^{\frac{q}{p}}}\nonumber\\
&\geq&\sup_{|a|>l}\frac{\mu_{u,\varphi}(\bigtriangleup(a,r))}{(1-|a|^2)^{\frac{q}{p}}}.\nonumber
\endr
Hence,
\begr
A+\epsilon&\geq&
\limsup_{|a|\to1}\frac{\mu_{u,\varphi}(\bigtriangleup(a,r))}{(1-|a|^2)^{\frac{q}{p}}}.\nonumber
\endr
Since $\epsilon$ is arbitrary, we obtain
\begr
\lim_{s\to1}\|\mu_{u,\varphi}|_{\D\backslash\D_s}\|_{p,q,r}^q&\geq&
\limsup_{|a|\to1}\frac{\mu_{u,\varphi}(\bigtriangleup(a,r))}{(1-|a|^2)^{\frac{q}{p}}}.\nonumber
\endr

 Now we consider the remainder of the proof.   Since $uC_\varphi: H^p\to L^q(\partial\D)$ is bounded, by Theorem 2.1, we have
$\mu_{u,\varphi}|_{\partial\D}=0$.  See \cite[Theorem 5]{CZ2} and the proof of \cite[Theorem 2]{CZ2} for the rest of the proof. Although in the proof it is assumed that  the function $u$ is analytic, the proof also work if it is only measurable.
 The comparability of the quantities follows from the proofs.\msk

 (ii) Let $\eta\in (0,1)$ be fixed. For $w\in\D$, let $K_w(z)=\frac{1}{1-\overline{w}z}$,~$z\in\D$. Then $K_w\in H^\infty\subset H^{p^\prime}$, where $1/p+1/p^\prime=1$.
 For $f\in H^p$ and $g\in H^{p^\prime}$, we denote
 $$\langle f,g\rangle=\int_{\partial\D}f(\xi)\overline{g(\xi)}dm(\xi).$$
It is easy to see that for every $f\in H^p$,
 $$f(w)=\langle f,K_w\rangle~~~\mbox{~~and~~} \langle R_nf,K_w\rangle = \langle f, R_nK_w\rangle.$$
 Thus,
 \begr
 |R_nf(w)| &=& |\langle R_nf,K_w\rangle| = |\langle f, R_nK_w\rangle|
 \leq \|f\|_p\|R_nK_w\|_\infty.\nonumber
 \endr
For all $\xi\in\varphi^{-1}(\D_\eta)$, let $w=\varphi(\xi)$. Then $|w|<\eta$.
Since
 $$R_nK_w(z) = R_n(\sum_{k=0}^\infty\overline{w}^kz^k) = \sum_{k=n+1}^\infty\overline{w}^kz^k,$$
one has
$$\|R_nK_w\|_\infty\leq\frac{\eta^{n+1}}{1-\eta}.$$
Therefore, 
$$\lim_{n\to\infty}\sup_{\|f\|_{p}\leq1}\int_{\varphi^{-1}(\D_\eta)}|(C_\varphi\circ R_nf)(\xi)|^qdm(\xi)\leq \lim_{n\to\infty}\frac{\eta^{n+1}}{1-\eta}=0.$$

\section{Proofs of main results }

To prove the main results in this paper, we need the following three lemmas.\msk

\noindent{\bf Lemma 3.1.}    {\it Let $0<r<1$. Then there exists a constant $C=C(r)>0$ such that whenever $a\in\D$ and $z\in\bigtriangleup(a,r)$,
\begr
\frac{|a|}{C}\rho(z,w)\leq \Big|1-\frac{1-\overline{a}z}{1-\overline{a}w}\Big|\leq C|a|\rho(z,w)\nonumber
\endr
for every $w\in \overline{\D}.$
}\msk

 {\it Proof.}  The proof is similarly with \cite[Lemma 4.3]{S1}. We only notice that
$|1-\overline{z}w|\approx|1-\overline{a}w|,$ whenever $a\in\D$, $z\in\bigtriangleup(a,r)$  and $w\in\overline{\D}$.\msk

\noindent{\bf Lemma 3.2.}    {\it Let $0<r<1$, $\gamma>0$. Then there exist constants $C_1=C_1(r,\gamma), C_2=C_2(r,\gamma)>0$ such that whenever $a\in\D$ and $z\in\bigtriangleup(a,r)$,
\begr
C_1\frac{|a|\rho(z,w)}{(1-|a|^2)^\gamma}\leq \Bigg|\bigg(\frac{1-|a|^2}{(1-\overline{a}z)^2}\bigg)^{\gamma}-\bigg(\frac{1-|a|^2}{(1-\overline{a}w)^2}\bigg)^{\gamma}\Bigg|\leq C_2\frac{|a|\rho(z,w)}{(1-|a|^2)^\gamma}\nonumber
\endr
for every $w\in \overline{\D}.$
}\msk

 {\it Proof.}  Let $a\in\D$, $w\in\overline{\D}$ and $z\in\bigtriangleup(a,r)$. By the proof of \cite[ Lemma 4.4]{S1}, we have
$$\bigg|1-\big(\frac{1-\overline{a}z}{1-\overline{a}w}\big)^{2\gamma}\bigg|\approx\bigg|1-\frac{1-\overline{a}z}{1-\overline{a}w}\bigg|.$$
Applying this, $|1-\overline{a}z|\approx 1-|a|^2$ and Lemma 2.1 we get the desired result. \msk

\noindent{\bf Lemma 3.3.} \cite{S1}  {\it Let $f\in H^1$ and $0<r<1$. Then there exist a constant $C=C(r)$ such that
\begr
|f(z)-f(a)|\leq C\rho(z,a)P|f|(a)\nonumber
\endr
for every $z\in \bigtriangleup(a,r).$ Here $Pf$ is the Poisson transformation of $f$, i.e.,
 $$Pf(z)=\frac{1}{2\pi}\int_0^{2\pi}\frac{1-|z|^2}{|1-\overline{z}e^{i\theta}|^2}f(e^{i\theta})d\theta,~~ f\in L^1(\partial\D).$$
}

Now we are in a position to prove our main results in this paper.\msk

{\it Proof of Theorem 1.1.} First, we consider the lower bound. Suppose that $C_\varphi-C_\psi: H^p\to H^q$ is bounded. Since $H^q$ is compact embedding into $H^p$, we have $C_\varphi-C_\psi:H^p\to H^p$ is compact.
Therefore, by \cite[Theorem 1 and Lemma 4]{NS}, $|\varphi|<1$ and $|\psi|<1$ a.e. on $\partial\D.$ Thus, $\mu_{\sigma,\varphi}(\partial\D)=0$.

Let $k_a(z)=\big(\frac{1-|a|^2}{(1-\overline{a}z)^2}\big)^{\frac{1}{p}}$. By Lemma 3.2, we get
\begr
&&\|C_\varphi-C_\psi\|_{H^p\to H^q}^q\nonumber\\
&\geq&\sup_{a\in\D}\|(C_\varphi-C_\psi)k_a\|_{q}^q\nonumber\\
&=&\sup_{a\in\D}\int_{\partial{\D}} \Bigg|\bigg(\frac{1-|a|^2}{(1-\overline{a}\varphi(\xi))^2}\bigg)^{\frac{1}{p}}-
\bigg(\frac{1-|a|^2}{(1-\overline{a}\psi(\xi))^2}\bigg)^{\frac{1}{p}}\Bigg|^qdm(\xi) \\
&\gtrsim&\sup_{|a|\geq 2^{-3}}\int_{\varphi^{-1}(\bigtriangleup(a,\frac{1}{2}))\cap\partial\D}\frac{|\sigma(\xi)|^q}{(1-|a|^2)^{\frac{q}{p}}}dm(\xi)\nonumber\\
&=&\sup_{|a|\geq 2^{-3}}\frac{\mu_{\sigma,\varphi}(\bigtriangleup(a,\frac{1}{2}))}{(1-|a|^2)^{\frac{q}{p}}}.\nonumber
\endr
Noting that $\bigtriangleup(a,2^{-3})\subset\bigtriangleup(2^{-2},2^{-1})$ for every $a\in \bigtriangleup(0,2^{-3})$, we get
 $$\sup_{|a|< 2^{-3}}\frac{\mu_{\sigma,\varphi}(\bigtriangleup(a,2^{-3}))}{(1-|a|^2)^{\frac{q}{p}}}\leq \frac{\mu_{\sigma,\varphi}(\bigtriangleup(2^{-2},\frac{1}{2}))}{(1-(2^{-2})^2)^{\frac{q}{p}}}.$$
Hence
$$\sup_{a\in\D}\frac{\mu_{\sigma,\varphi}(\bigtriangleup(a,2^{-3}))}{(1-|a|^2)^{\frac{q}{p}}}\lesssim\sup_{a\in\D}\|(C_\varphi-C_\psi)k_a\|_{q}^q.$$
Thus, by Theorem 2.1, we obtain that
$\sigma C_\varphi:H^p\to L^q(\partial\D)$ is bounded and
$$\|C_\varphi-C_\psi\|_{H^p \rightarrow H^q}^q \gtrsim\sup_{a\in\D}\|(C_\varphi-C_\psi)k_a\|_{q}^q\gtrsim \|\sigma C_\varphi\|_{H^p\to L^q(\partial\D)}^q .$$
Similarly,
$$\|C_\varphi-C_\psi\|_{H^p \rightarrow H^q}^q \gtrsim\sup_{a\in\D}\|(C_\varphi-C_\psi)k_a\|_{q}^q\gtrsim \|\sigma C_\psi\|_{H^p\to L^q(\partial\D)}^q,$$
and hence
 \begr \|C_\varphi-C_\psi\|_{H^p \rightarrow H^q}^q \gtrsim  \|\sigma C_\varphi\|_{H^p\to L^q(\partial\D)}^q + \|\sigma C_\psi\|_{H^p\to L^q(\partial\D)}^q.
 \endr

Next we consider the upper bound.
\begr
&&\|C_\varphi-C_\psi\|_{H^p\to H^q}^q\nonumber\\
 &=&\sup_{\|f\|_p\leq 1}\|(C_\varphi-C_\psi)f\|_{q}^q\nonumber\\
&=&\sup_{\|f\|_p\leq 1}\bigg(\int_{|\sigma(\xi)|\geq \frac{1}{2}}+\int_{|\sigma(\xi)|<\frac{1}{2}}\bigg)|f\circ\varphi(\xi)-f\circ\psi(\xi)|^qdm(\xi)\nonumber\\
&:=&I_1+I_2.\nonumber
\endr
It is easy to see that
\begr
I_1 \lesssim \|\sigma C_\varphi\|_{H^p\to L^q(\partial\D)}^q+\|\sigma C_\psi\|_{H^p\to L^q(\partial\D)}^q.\nonumber
\endr
By Lemma 2.3,
\begr
I_2
&=&\sup_{\|f\|_p\leq 1}\int_{|\sigma(\xi)|<\frac{1}{2}}|f\circ\varphi(\xi)-f\circ\psi(\xi)|^qdm(\xi)\nonumber\\
&\leq&\sup_{\|f\|_p\leq 1}\int_{|\sigma(\xi)|<\frac{1}{2}}|\sigma(\xi)|^q\big(P|f|\circ\varphi(\xi)\big)^qdm(\xi)\nonumber\\
&\leq&\sup_{\|f\|_p\leq 1}\int_{\overline{\D}}\big(P|f|(z)\big)^qd\mu_{\sigma,\varphi}(z). \nonumber
\endr
Let $\tilde{g}(z)$ denote the harmonic conjugate function of $g(z):=P|f|(z)$, normalized so that $\tilde{g}(0)=0$, and let $v=g+i\tilde{g}$. Then both $g$ and $\tilde{g}$ are belong to the harmonic Hardy space $h^p$ and hence $v\in H^p$. By M. Riesz theorem (see \cite[Theorem 2.3 of Chapter 3 ]{Ga} or \cite[Theorem 4.1]{Duren}) and Minkowski inequality,  we get $$\|v\|_p\approx\|g\|_p+\|\tilde{g}\|_p\lesssim \|g\|_p\leq \|f\|_p.$$
 By Theorem 2.1,
\begr \int_{\overline{\D}}\big(P|f|(z)\big)^qd\mu_{\sigma,\varphi}(z)
&\leq&\int_{\overline{\D}}|v(z)|^qd\mu_{\sigma,\varphi}(z)\nonumber\\
&\lesssim& \|\sigma C_\varphi\|_{H^p\to L^q(\partial\D)}^q\|v\|_{p}^q\nonumber\\
&\lesssim&
\|\sigma C_\varphi\|_{H^p\to L^q(\partial\D)}^q \|f\|_{p}^q.\nonumber \endr
Thus, $$I_2\lesssim\|\sigma C_\varphi\|_{H^p\to L^q(\partial\D)}^q .$$
Hence
\begr \|C_\varphi-C_\psi\|_{H^p\to H^q}^q  \lesssim \|\sigma C_\varphi\|_{H^p\to L^q(\partial\D)}^q+ \|\sigma C_\psi\|_{H^p\to L^q(\partial\D)}^q.
\endr
 Therefore,  by (1), (2) and (3), we get
\begr \|C_\varphi-C_\psi\|_{H^p\to H^q}^q& \approx&\sup_{a\in\D}\|(C_\varphi-C_\psi)k_a\|_{q}^q\nonumber\\
&\approx &\|\sigma C_\varphi\|_{H^p\to L^q(\partial\D)}^q+ \|\sigma C_\psi\|_{H^p\to L^q(\partial\D)}^q.\nonumber
\endr
The result follows by Theorem 2.1.\msk

{\it Proof of Theorem 1.2.} We can assume that $\varphi\neq\psi$.
Let $\{a_n\}$ be any sequence in $\D$ such that $a_n\to 1$ as $n\to \infty$.
Let $k_{a_n}(z)=\big(\frac{1-|a_n|^2}{(1-\overline{a_n}z)^2}\big)^{\frac{1}{p}}$. Then $\|k_{a_n}\|_{H^p}=1$ and $k_{a_n}\rightarrow0$ weakly in $H^p$ as $n\to\infty$. Let $S$ be a compact operator from $H^p$ into $H^q$. Then $\lim_{n\rightarrow\infty}\|Sk_{a_n}\|_{H^q}=0$. Hence,
\begr
\|C_{\varphi}-C_{\psi}-S\|_{H^p\to H^q} &\geq& \limsup_{n\rightarrow\infty}\|(C_{\varphi}-C_{\psi}-S)k_{a_n}\|_{q}\nonumber\\
&\geq& \limsup_{n\rightarrow\infty}\|(C_{\varphi}-C_{\psi})k_{a_n}\|_{q}.
\nonumber
\endr
Thus,
 \begr
\|C_{\varphi}-C_{\psi}\|_{e,H^p\to H^q }\geq\limsup_{|a|\rightarrow1}\|(C_{\varphi}-C_{\psi})k_{a}\|_{q}.
\endr
Now, we prove that
 \begr
\limsup_{|a|\rightarrow1}\|(C_{\varphi}-C_{\psi})k_{a}\|_{q}^q\gtrsim \|\sigma C_{\varphi}\|_{e, H^p\to L^q(\partial\D)}^q+\|\sigma C_{\psi}\|_{e,H^p\to L^q(\partial\D)}^q.
\endr
 Since $C_\varphi-C_\psi: H^p\to H^q$ is bounded, by Theorem 1.1, we have $\sigma C_\varphi:H^p\to L^q(\partial\D)$ is bounded.
Let $0<r<1$. By Lemma 3.2 and Theorem 2.2, we obtain
\begr
&&\limsup_{|a|\rightarrow1}\|(C_{\varphi}-C_{\psi})k_{a}\|_{q}^q\nonumber\\
&=&\limsup_{|a|\to 1}\int_{\partial{\D}} \Bigg|\bigg(\frac{1-|a|^2}{(1-\overline{a}\varphi(\xi))^2}\bigg)^{\frac{1}{p}}-
\bigg(\frac{1-|a|^2}{(1-\overline{a}\psi(\xi))^2}\bigg)^{\frac{1}{p}}\Bigg|^qdm(\xi)\nonumber\\
&\gtrsim&\limsup_{|a|\rightarrow1}\int_{\varphi^{-1}(\bigtriangleup(a,r))\cap\partial\D}
\frac{|\sigma(\xi)|^q}{(1-|a|^2)^{\frac{q}{p}}}dm(\xi)\nonumber\\
&=&\limsup_{|a|\rightarrow1}\frac{\mu_{\sigma,\varphi}(\bigtriangleup(a,r))}{(1-|a|^2)^{\frac{q}{p}}}\nonumber\\
&\approx& \|\sigma C_{\varphi}\|_{e, H^p\to L^q(\partial\D)}^q.
\nonumber
\endr
Similarly, we get
\begr
\limsup_{|a|\rightarrow1}\|(C_{\varphi}-C_{\psi})k_{a}\|_{q}^q\gtrsim \|\sigma C_{\psi}\|_{e,H^p\to L^q(\partial\D)}^q.
\nonumber
\endr
Finally, we prove that
$$
\|C_{\varphi}-C_{\psi}\|_{e,H^p\to H^q }^q\lesssim \|\sigma C_{\varphi}\|_{e, H^p\to L^q(\partial\D)}^q+\|\sigma C_{\psi}\|_{e,H^p\to L^q(\partial\D)}^q.
$$
Since the partial sum operator $S_n$ is compact,  we get
$$\|C_{\varphi}-C_{\psi}\|_{e, H^p\to H^q}\leq \limsup_{n\to\infty}\|(C_\varphi-C_\psi)R_n\|_{H^p\to L^q(\partial\D)}.
 $$
Denote $E=\{\xi\in\partial\D:|\sigma(\xi)|\geq 1/2\}$ and $E^\prime=\partial\D\backslash E$. Then
\begr
I_n(f)&:=&\int_{E}|(C_\varphi-C_\psi)R_nf(\xi)|^qdm(\xi)\nonumber\\
&\leq&2^q\Big(\int_{E}|(\sigma C_\varphi)R_nf(\xi)|^qdm(\xi)+\int_{E}|(\sigma C_\psi)R_nf(\xi)|^qdm(\xi)\Big)\nonumber\\
&\leq&2^q\Big(\|(\sigma C_\varphi)R_n\|_{H^p\to L^q(\partial\D)}^q+\|(\sigma C_\psi)R_n\|_{H^p\to L^q(\partial\D)}^q\Big),\nonumber
\endr
whenever $\|f\|_p\leq 1$ and $n\in\N$. Thus by Theorem 2.2 $(i)$,
\begr
\limsup_{n\to\infty}\sup_{\|f\|_p\leq 1}I_n(f)\lesssim \|\sigma C_{\varphi}\|_{e,H^p\to L^q(\partial\D)}^q+\|\sigma C_{\psi}\|_{e,H^p\to L^q(\partial\D)}^q.\nonumber
\endr
Denote
$$J_n(f):=\int_{E^\prime}|(C_\varphi-C_\psi)R_nf(\xi)|^qdm(\xi).$$
For all $a, z, w\in\D$, from \cite[Lemma 1.4 of Chapter 1]{Ga} or \cite{sh} we see that
$$\rho(z,w)\leq \frac{\rho(z,a)+\rho(a,w)}{1+\rho(z,a)\rho(a,w)}.  $$
Let $s\in (0,1)$ be arbitrary. Suppose $z\in E^\prime\cap\varphi^{-1}(\D_s)$. By the last inequality we can find $s^\prime=\frac{\frac{1}{2}+r}{1+\frac{r}{2}}\in(0,1)$ such that $E^\prime\cap\varphi^{-1}(\D_s)\subset\psi^{-1}(\D_{s^\prime})$. Thus by Theorem 2.2 $(ii)$,
 $$\lim_{n\to\infty}\sup_{\|f\|_p\leq1}\int_{E^\prime\cap\varphi^{-1}(\D_s)}|(C_{\varphi}\circ R_nf)(\xi)|^q dm(\xi)=0$$
 and
  $$\lim_{n\to\infty}\sup_{\|f\|_p\leq1}\int_{E^\prime\cap\varphi^{-1}(\D_s)}|(C_\psi\circ R_nf)(\xi)|^q dm(\xi)=0.$$
Hence
\begr
\limsup_{n\to\infty}\sup_{\|f\|_p\leq1}J_n(f)&\lesssim&\limsup_{n\to\infty}\sup_{\|f\|_p\leq1}\int_{F}|(C_\varphi-C_\psi)\circ R_nf(\xi)|^q dm(\xi)\nonumber\\
&\lesssim&\sup_{\|f\|_p\leq1}\int_{F}|(C_\varphi-C_\psi) f(\xi)|^q dm(\xi),\nonumber
\endr
 where $F=E^\prime\cap\varphi^{-1}(\overline\D\backslash\D_s)$ and we used the fact that the operators $S_n$ are uniformly bounded (see \cite[Proposition 1]{Z}), so does $R_n$.

 Using Lemma 3.3, we get
  \begr
 \int_F|(C_\varphi-C_\psi)f(\xi)|^qdm(\xi)& \lesssim&\int_{F}|\sigma(\xi)|^q\big(P|f|\circ\varphi(\xi)\big)^q
dm(\xi)\nonumber\\
&=&\int_{\varphi(F)}\big(P|f|(z)\big)^q
d\mu_{\sigma,\varphi}(z)\nonumber\\
&\leq&\int_{\overline\D\backslash\D_s}\big(P|f|(z)\big)^q
d\mu_{\sigma,\varphi}(z).\nonumber
\endr
Let $\tilde{g}(z)$ denote the harmonic conjugate function of $g(z):=P|f|(z)$ with $\tilde{g}(0)=0$ and let $v=g+i\tilde{g}$. Then $v\in H^p$, $|g(z)|\leq |v(z)|$, and $$\|v\|_p\approx\|g\|_p\leq\|f\|_p\leq 1.$$
Therefore,
  \begr
 \int_{\overline\D\backslash\D_s}\big(P|f|(z)\big)^q d\mu_{\sigma,\varphi}(z)
\lesssim\int_{\overline\D\backslash\D_s}|v(z)|^qd\mu_{\sigma,\varphi}(z).
\nonumber
 \endr
 Since $C_\varphi-C_\psi:H^p\to H^q$ is bounded, we have $\sigma C_\varphi:H^p\to L^q(\partial\D)$ is bounded. Thus,
 $\mu_{\sigma,\varphi}|_{\partial\D}=0$ and
$$ \|\big(\mu_{\sigma,\varphi}|_{\overline\D\backslash\D_s}\big)|_{\D}\|_{p,q,r}\leq \|\mu_{\sigma,\varphi}\|_{p,q,r}<\infty. $$
This show that $\mu_{\sigma,\varphi}|_{\overline\D\backslash\D_s}$ is a $\frac{q}{p}$-Carleson measure on $\overline\D$. Then,
 \begr
 \int_F|(C_\varphi-C_\psi)f(\xi)|^qdm(\xi) &\lesssim&\int_{\overline\D\backslash\D_s}|v(z)|^q
d\mu_{\sigma,\varphi}(z)\nonumber\\
&=&\int_{\overline\D}|v(z)|^qd(\mu_{\sigma,\varphi}|_{\overline\D\backslash\D_s})(z)\nonumber\\
&\lesssim&\|(\mu_{\sigma,\varphi}|_{\overline\D\backslash\D_s})|_{\D}\|^q_{p,q,r}\nonumber\\
&=&\|\mu_{\sigma,\varphi}|_{\D\backslash\D_s}\|^q_{p,q,r}.\nonumber
 \endr
Letting $s\to1$, we get
 $$\limsup_{n\to\infty}\sup_{\|f\|_p\leq1}J_n(f)\lesssim\|\sigma C_\varphi\|_{e, H^p\to L^q(\partial\D)}^q.$$
 Therefore,
\begr
\|C_{\varphi}-C_{\psi}\|_{e, H^p\to H^q}^q&\leq& \limsup_{n\to\infty}\sup_{\|f\|_p\leq1}\|(C_\varphi-C_\psi)R_n\|_{H^p\to L^q(\partial\D)}^q\nonumber\\
&=&\limsup_{n\to\infty}\sup_{\|f\|_p\leq1}I_n(f)+\limsup_{n\to\infty}\sup_{\|f\|_p\leq1}J_n(f)\nonumber\\
&\lesssim&\|\sigma C_\varphi\|_{e,H^p\to L^q(\partial\D)}^q+\lesssim\|\sigma C_\psi\|_{e, H^p\to L^q(\partial\D)}^q.
\endr
By (4), (5), (6)   and Theorem 2.2, we get the desired result. \msk

\noindent{\bf Remark 3.1.} { Theorem 1.2 motivate us to study a possible connection between the differences of composition operators and the corresponding weighted composition operators on $H^p$. We conjecture that
\begr
&&\|C_\varphi-C_\psi\|_{e, H^p\to H^p}^p\nonumber\\
&\approx&\limsup_{|a|\to 1}\int_{\partial{\D}} \Big|\bigg(\frac{1-|a|^2}{(1-\overline{a}\varphi(\xi))^2}\bigg)^{\frac{1}{p}}-
\bigg(\frac{1-|a|^2}{(1-\overline{a}\psi(\xi))^2}\bigg)^{\frac{1}{p}}\Big|^pdm(\xi)\nonumber\\
&\approx&\limsup_{|a|\to 1}\int_{\partial{\D}} \Big(\frac{1-|a|^2}{|1-\overline{a}\varphi(\xi)|^2}+\frac{1-|a|^2}{|1-\overline{a}\psi(\xi)|^2}\Big)
\rho(\varphi(\xi),\psi(\xi))^pdm(\xi).\nonumber\endr
}

\noindent{\bf Remark 3.2.} In \cite{NS}, Nieminen and  Saksman showed that the compactness of $C_\varphi-C_\psi$ on $H^p$ is independent of $p\in[1,\infty)$. Therefore, we also conjecture that $C_\varphi-C_\psi:H^p\to H^p$$(1\leq p<\infty)$ is compact
if and only if
\begr
\lim_{|a|\to 1}\int_{\partial{\D}} \Big(\frac{1-|a|^2}{|1-\overline{a}\varphi(\xi)|^2}+\frac{1-|a|^2}{|1-\overline{a}\psi(\xi)|^2}\Big)
\rho(\varphi(\xi),\psi(\xi))dm(\xi)=0.\nonumber\endr\msk

\noindent{\bf Remark 3.3.} In \cite{Sj1}, Shapiro showed that the  operator $C_\varphi$ is compact on $H^p$ if and only if
$$\lim_{|z|\to1}\frac{N_{\varphi}(z)}{\log{\frac{1}{|z|}}}=0,$$
where $N_{\varphi}$, called the Nevanlinna counting function,  is defined by $$N_\varphi(z)=\sum_{a\in \varphi^{-1}(z)}\log\frac{1}{|a|},  ~~~z\in \D\setminus\varphi(0).$$
We ask whether $C_\varphi-C_\psi:H^p\to H^p$ is compact if and only if
$$\lim_{|z|\to1}\rho(\varphi(z),\psi(z))\bigg(\frac{N_\varphi(z)}{\log{\frac{1}{|z|}}}+\frac{N_\psi(z)}{\log{\frac{1}{|z|}}}\bigg)=0?$$

{\bf Acknowledgement.}  This  project was partially supported by NSF of China (No.11471143 and No.11720101003). \msk

\end{document}